\documentclass[12pt,reqno]{amsart}
\usepackage{amsmath}
\usepackage{amssymb}
\usepackage{amsthm}
\usepackage[left=2.9cm,top=2.5cm,right=2.9cm,bottom=2.6cm]{geometry}
\usepackage{graphicx}
\usepackage[utf8]{inputenc}
\usepackage[T1]{fontenc}
\usepackage{url}

\begin{document}
\newtheorem{theorem}{Theorem}[section]
\newtheorem{lemma}[theorem]{Lemma}
\newtheorem{definition}[theorem]{Definition}
\newtheorem{conjecture}[theorem]{Conjecture}
\newtheorem{proposition}[theorem]{Proposition}
\newtheorem{algorithm}[theorem]{Algorithm}
\newtheorem{corollary}[theorem]{Corollary}
\newtheorem{observation}[theorem]{Observation}
\newtheorem{remark}[theorem]{Remark}
\newtheorem{problem}[theorem]{Problem}
\newtheorem{example}[theorem]{Example}
\newcommand{\noin}{\noindent}
\newcommand{\ind}{\indent}
\newcommand{\al}{\alpha}
\newcommand{\om}{\omega}
\newcommand{\pp}{\mathcal P}
\newcommand{\ppp}{\mathfrak P}
\newcommand{\R}{{\mathbb R}}
\newcommand{\N}{{\mathbb N}}
\newcommand{\Z}{{\mathbb Z}}
\newcommand\eps{\varepsilon}
\newcommand{\E}{\mathbb E}
\newcommand{\Prob}{\mathbb{P}}
\newcommand{\pl}{\textrm{C}}
\newcommand{\dang}{\textrm{dang}}
\renewcommand{\labelenumi}{(\roman{enumi})}
\newcommand{\bc}{\bar c}
\newcommand{\G}{{\mathfrak S}}
\newcommand{\T}{{\mathfrak T}}
\newcommand{\dist}{\textrm{dist}}
\newcommand{\rad}{\textrm{rad}}
\newcommand{\Remark}[1]{\par \noindent \textnormal{\textbf{Note.~}} #1}

\title{The 2-surviving rate of planar graphs with average degree lower than $\frac{9}{2}$}
\author{Przemys\l{}aw Gordinowicz}
\address{Institute of Mathematics, Lodz University of Technology, \L{}\'od\'z, Poland}
\email{pgordin@p.lodz.pl}
\begin{abstract}
Let $G$ be any connected graph on $n$ vertices, $n \ge 2.$ Let $k$ be any positive integer. Suppose that a fire breaks out at some vertex of $G.$ Then, in each turn firefighters can protect at most $k$ vertices of $G$ not yet on fire; Next the fire spreads to all unprotected neighbours of burning vertices. 
The \emph{$k$-surviving rate} of G, denoted by  $\rho_k(G),$ is the expected fraction of vertices that can be saved from the fire, provided that the starting vertex is chosen uniformly at random. 

In this note, it is shown that for any planar graph $G$ with average degree $\frac{9}{2} - \eps,$ where $\eps \in (0, 1],$ there is $\rho_2(G) \ge \frac{2}{9}\eps.$ 
In particular, the result implies a significant improvement of the bound for 2-surviving rate for triangle-free planar graphs (Esperet, van den Heuvel, Maffray and Sipma~\cite{esperet}) and for planar graphs without 4-cycles (Kong, Wang, Zhang~\cite{kwzhang}). The proof is done using the separator theorem for planar graphs.

This paper is the corrected version of~\cite{pg_planarfire2} unified with the corrigendum~\cite{pg_planarfire2_err}. 
\end{abstract}

%\linenumbers

\vspace{150pt}

\maketitle

\section{Introduction\label{sec:intro}}
The following \emph{Firefighter Problem} was introduced by Hartnell~\cite{hartnell}. Let $G$ be any connected graph on $n$ vertices, $n \ge 2.$ Let $k$ be any positive integer. Suppose that a~fire breaks out at some vertex $v \in V(G).$ Then in each turn firefighters can protect $k$ vertices of $G,$ not yet on fire and the protection is permanent. Next the fire spreads to all the unprotected neighbours of vertices already on fire. The process ends when the fire can no longer spread. The goal is to save as much as possible and the question is how many vertices can be saved. We refer the reader to the survey by Finbow and MacGillivray~\cite{fm} for more information on the background of the problem and directions for its consideration. 

In this note we focus on the following aspect of the problem. 
Let $sn_k(G, v)$ denote the maximum number of vertices of $G$ that $k$ firefighters can save when the fire breaks out at the vertex $v.$ This parameter may depend heavily on the choice of the starting vertex $v,$ for example when the graph $G$ is a star. Therefore Cai and Wang~\cite{cai} introduced the following graphical parameter: the \emph{$k$-surviving rate} $\rho_{k}(G)$ is the expected fraction of vertices that can be saved by $k$ firefighters, provided that the starting vertex is chosen uniformly at random, that is $$\rho_k(G) = \frac{1}{|V(G)|^2} \sum_{v \in V(G)} sn_k(G, v).$$
Given class of graphs the goal is to bound the surviving rate from below. We note that an equivalent graph parameter, namely the expected damage, $ed(G) = |V(G)|(1-\rho_1(G))$ was considered earlier by Finbow, Hartnell, Li and Schmeisser~\cite{fhls}. Oppositely, for several graph classes they investigated representants which minimises the expected damage.  

It is not surprising that the surviving rate is connected with the density of a~graph. Prałat~\cite{pralat3011,pralatfire} has provided a~threshold for the average degree, which guarantees a~positive surviving rate with a given number of firefighters. Precisely, for $k \in \N^+$ let us define
$$\tau_k = \left\{ \begin{array}{lcl}
\frac{30}{11} & \textrm{~for~} & k = 1\\
k + 2 - \frac{1}{k+2}& \textrm{~for~} & k \ge 2.\\    
\end{array}\right.$$
Then there exists a constant $c > 0,$ such that for any $\epsilon > 0,$ any $n \in \N^+$ and any graph $G$ on $n$ vertices and at most $(\tau_k - \epsilon)n/2$ edges one has $\rho_k(G) > c \cdot \epsilon >0.$ Moreover, there exists a~family of graphs with the average degree tending to $\tau_k$ and the $k$-surviving rate tending to 0, which shows that the above result is the best possible. In particular, Pra\l{}at's results yield that graphs with the average degree lower than $\frac{15}{4} - \eps$ have the 2-surviving rate at least $\frac{8}{75}\eps.$

While settled in general case, the $k$-surviving rate is still being investigated for particular families of graphs --- the most important is the case of planar graphs. Cai and Wang~\cite{cai} asked about the minimum number $k$ such that $\rho_k(G) > c$ for some positive constant $c$ and any planar graph $G.$ It is easy to see that $\rho_1(K_{2,n}) \overrightarrow{_{\ n \to \infty}\ } \, 0,$ hence $k \ge 2.$ So far, the best known upper bound for this number  is 3: 
\begin{theorem}[\cite{pg_planarfire}] \label{thm:3f}
Let $G$ be any planar graph. Then
$$\rho_{3}(G) > \frac{2}{21}.$$
\end{theorem}

It has been shown that 2 is the upper bound for triangle-free planar graphs~\cite{esperet} and planar graphs without 4-, 5- or 6-cycles~\cite{kwzhang, kwwu, fww}, respectively. In particular Esperet, van den Heuvel, Maffray and Sipma~\cite{esperet} have proven that
$$\rho_2(G) > \frac{1}{723636}$$ for any non-trivial triangle-free planar graph, while Kong, Wang and Zhang~\cite{kwzhang} that
$$\rho_2(G) > \frac{1}{76}$$ for any non-trivial planar graph without 4-cycles.

In this paper, we improve the above bounds with the following theorem:
\begin{theorem} \label{thm:main}
Let $G$ be any connected planar graph with $n \ge 2$ vertices and $m$ edges. If for some $\epsilon \in (0, \frac{5}{2}]$ one has $\frac{2m}{n} = \frac{9}{2} - \epsilon$ then
$$ \rho_{2}(G) \ge \left\{ \begin{array}{ll}
\frac{2}{9}\eps & \textrm{~for~} \epsilon \le 1\\
\\
\frac{2}{9}\eps - \frac{1}{n} & 	\textrm{~otherwise.} \\
\end{array}\right.  $$
\end{theorem}

In particular, because triangle-free planar graphs have average degree lower than 4, while planar graphs without 4-cycles have it bounded by $\frac{30}{7}$ (see eg.~\cite[Lemma 1.7]{mop} or derive it directly from Euler Formula), we get immediately
\begin{corollary} \label{cor:trianglefree}
Let $G$ be any triangle-free planar graph on at least 2 vertices. Then
$$\rho_{2}(G) > \frac{1}{9}.$$
\end{corollary}
\begin{corollary} \label{cor:4cycles}
Let $G$ be any graph on at least 2 vertices. If $G$ does not contain any 4-cycle then $$\rho_{2}(G) > \frac{1}{21}.$$
\end{corollary}

The paper is organised as follows. Section~\ref{sec:separators} contains a brief description of the firefighters' strategy developed from the separator theorem of planar graphs. Section~\ref{sec:proof} contains the proof of Theorem~\ref{thm:main}. 
For background on graph theory we refer the reader to the book \cite{diestel}, in particular to Chapter 4 for the terminology and basics regarding planarity and plane embeddings of graphs. Note that according to~\cite{diestel}, we imagine the edges of a~plane graph as simple polygonal curves, while by a Jordan curve we mean a simple, closed polygonal curve. Considering particular plane graph $G$ and some Jordan curve $C$ by $\textrm{in} C,$ $\overline{\textrm{in} C},$ $\textrm{ex} C$ and $\overline{\textrm{ex} C}$ we denote sets of vertices of $G$ lying in the interior region of $C,$ in the closed interior region (including $C$), in the exterior region and in the closed exterior region of $C,$ respectively. 

\section{Separators and the firefighters' strategy\label{sec:separators}}

The proof is done using the lemma given by Lipton and Tarjan to prove the separator theorem for planar graphs~\cite{liptontarjan}. The key lemma in their proof, reformulated for our purpose, is quoted below. A similar approach --- using the lemma of Lipton and Tarjan to the firefighter problem on planar graphs, was first applied by Floderus, Lingas and Persson~\cite{lingas2}, with a slightly different notation of approximation algorithms. They have proven a theorem analogous to Lemma~\ref{lem:many}. This method was also used to prove Theorem~\ref{thm:3f}. The proof and some discussion on such an approach is included in~\cite{pg_planarfire}.

\begin{lemma} \label{thm:separator}
Let $G$ be any connected $n$-vertex plane graph $(n \ge 2)$ and $T$ be any spanning tree of $G.$ Then there exists an arc $\widetilde{uv}$ between two vertices $u, v \in V(G),$ not crossing any edge from $E(G),$ such that the unique Jordan curve $C,$ consisting of $\widetilde{uv}$ and some edges of $T,$ has the property that the number of vertices in the interior region of $C$ as well as in the exterior region of $C$ is lower than $\frac{2}{3}n.$
\end{lemma}
Actually in~\cite{liptontarjan} there is considered any spanning supergraph of $G$, say $H$, being a plane triangulation and the original lemma guarantees existence of an edge $uv \in E(H)$ with the desired property. Hence we may view an arc $\widetilde{uv}$ as some edge of some planar supergraph of $G$, while $uv$ may be an edge of $G$ or not. To construct the curve $C$, for any vertex $r \in V(G)$, let us consider two $T$-paths from $r$ to $u$ and $v,$ respectively. Let $z$ be the last common vertex on these paths. Then the curve  $C$ consists of the $T$-path from $z$ to $u,$ the arc $\widetilde{uv}$ and the $T$-path from $v$ to $z.$ 

The defending strategy is the following. Let $G$ be any $n$-vertex connected plane graph, where $n \ge 5.$ Suppose that the fire breaks out at a vertex $r.$ Consider a tree $T$ obtained by the breadth-first-search algorithm starting from the vertex $r.$ By Lemma~\ref{thm:separator} there is an arc $\widetilde{uv}$ joining vertices $u$ and $v$ determining the Jordan curve $C$ consisting of $\widetilde{uv}$ and some edges of $T$ such that $|\overline{\textrm{in} C}| > \frac{1}{3}n$ and $|\overline{\textrm{ex} C}| > \frac{1}{3}n.$ Note that the curve $C$ contains at most 2 vertices at any given distance from $r.$ This holds true because the curve $C$ is constructed from two shortest paths (as $T$ is a breadth-first-search tree), say $u_0u_1\dots u_i$ and $v_0v_1\dots v_j,$ where $v_0 = u_0 = r,$ $u_i = u$ and $v_j = v.$ The firefighters' strategy relies on protecting in the $t$th round, for $t \ge 1,$ vertices $u_t$ and $v_t,$ until every vertex from both paths, except $r,$ is protected. When the vertex $r$ does not belong to the curve $C$ then firefighters, protecting the whole separator, can save all the vertices in either $\overline{\textrm{in} C}$ or $\overline{\textrm{ex} C}.$ When the vertex $r$ belongs to the curve $C$ the protection along the separator may be not enough, as the fire may spread through the neighbours of $r$ located in the interior as well as in the exterior region of the curve $C.$
Because then either $\textrm{in} C$ or $\textrm{ex} C$ contains not more than $\left\lfloor\frac{\deg r -2}{2}\right\rfloor$ neighbours of $r,$ we get immediately Lemma~\ref{lem:many}.

To increase clarity of the presentation we decided to describe the firefighters' strategy in a way that is not strictly precise. It may happen that the vertex that should be protected according to our description does not exist or is already protected. 
For example, in the above paragraph, when $u_1 = v_1$ or in the $i$th round when $i > j.$ In such cases firefighters may either protect a random vertex instead or skip the protection at all.

\begin{lemma} \label{lem:many}
Let $G$ be any $n$-vertex plane graph, where $n \ge 5.$ Suppose that the fire breaks out at some vertex $r.$ Then using $2 + \left\lfloor\frac{\deg r -2}{2}\right\rfloor$ firefighters at the first step and 2 at the subsequent steps one can save more than $\frac{n}{3} - 1$ vertices. 
\end{lemma}

Note that, due to a divisibility condition, ``more than $\frac{n}{3}-1$'' means ``at least $\lceil \frac{n-2}{3} \rceil$''.
For any planar graph with average degree lower than 4 (e.g. triangle-free planar graphs) we get:

\begin{corollary} \label{cor:f32}
Let $G$ be any planar graph with average degree lower than 4, having at least 2 vertices of degree lower than 4. There is $$\rho_{3,2}(G) > \frac{1}{3}.$$
\end{corollary}
By $\rho_{3,2}(G)$ we denote the surviving rate of a graph $G$ in a scenario with 3 vertices protected in the first round and 2 in subsequent rounds. 
\begin{proof}
Let $G$ be any $n$-vertex planar graph with average degree lower than 4, having $x$ vertices of degree lower than 4 and $y$ vertices of degree at least 6. We have $x \ge 2$ and $y < \frac{3x}{2}$. Note that the first round is enough to save $n-1$ vertices when fire breaks out at vertex of degree lower than 4 and to save 3 vertices when degree of starting vertex is at least 6. The first condition implies that $\rho_{3,2}(G) > \frac{1}{3}$  for $n \le 4$, so we assume further $n > 4$. When fire breaks out at vertex of degree 4 or 5 by Lemma~\ref{lem:many} it is possible to save at least $\lceil \frac{n-2}{3} \rceil$ vertices. Altogether, the surviving rate $\rho_{3,2}(G)$ can be bounded from below by the solution $\alpha_0$ of the following integer program (for $x, y \in \N$):   
\begin{eqnarray} \label{eq:lp2}
\nonumber \textrm{minimize} & & \alpha = \frac{1}{n^2} \Big((n-1)x + 3y + \frac{n-2}{3}(n - x - y) \Big)\\
 \textrm{with conditions} & & y < \frac{3x}{2},\\
\nonumber & & x + y \le n \\
\nonumber & & x \ge 2.
\end{eqnarray}
One can calculate that the minimum is obtained for $x = y = 2$ and then $$\rho_{3,2}(G) \ge \alpha_0 = \frac{1}{3} +  \frac{16}{3n^2} > \frac{1}{3}.$$
\end{proof}

\section{The proof\label{sec:proof}}

Let us start the proof of Theorem~\ref{thm:main} with a simple observation derived from Lemma~\ref{lem:many}.
\begin{observation} \label{obs:simple}
Let $G$ be any $n$-vertex plane graph, where $n \ge 5.$ Let $r \in V(G)$ be a~vertex of degree at most $3.$ Then 
$$sn_{2}(G, r) > \left\{ \begin{array}{lcl}
n-1.01 & \textrm{~if~} & \deg(r) \le 2\\
n/3-1 & \textrm{~if~} & \deg(r) = 3\\
1 & \textrm{~if~} & \deg(r) > 4\\  
\end{array}\right.$$
\end{observation}

Dealing with vertices of degree 4 is a bit more complicated. We prove that there is a~defending strategy to save a large part of a graph unless more than one neighbour of the vertex where the fire broke out has a high degree.  
\begin{lemma} \label{lem:2f}
Let $G$ be any $n$-vertex planar graph, where $n \ge 5.$ Let $r \in V(G)$ be a vertex of degree $4$ with at most one neighbour of degree higher than 5. Then  
$$sn_{2}(G, r) > \frac{1}{6}n - 1.$$
\end{lemma}

\begin{proof}
Let $G$ be any $n$-vertex planar graph, where $n \ge 5,$ and let plane graph $H$ be its plane embedding. Let $r \in V(H)$ be a vertex of degree $4.$ Consider a~tree $T$ obtained by the breadth-first-search algorithm starting from the vertex $r.$ By Lemma~\ref{thm:separator} there are vertices $u$ and $v$ and the Jordan curve $C \subseteq T + \widetilde{uv}$ such that $|\overline{ \textrm{in} C}| > \frac{1}{3}n$ and $|\overline{\textrm{ex} C}| > \frac{1}{3}n.$ 
Let $u_0u_1\dots u_i$ and $v_0v_1\dots v_j$ be $T$-paths from $r = u_0 = v_0$ to $u = u_i$ and $v = v_j,$ respectively.

Consider an open disk around $r$ small enough to contain only one line segment per each of 4 edges adjacent to $r.$ Enumerate neighbours of $r$ as $a, b, c, d$ using a clockwise ordering applied to the corresponding line segments, so that $a = u_1.$ When $v_1 = a,$ using the strategy described as the introduction to Lemma~\ref{lem:many}, two firefighters can save more than $\frac{1}{3}n+1$ vertices of a graph (one extra in the first round as $u_1 = v_1$), because vertex $r$ does not belong to the curve $C.$ Similarly, when $v_1 \in \{b,d\},$ it is possible to save more than $\frac{1}{3}n-1$ vertices.  

Therefore suppose now that $v_1 = c.$ If in the graph $H$ there exists some shortest path from vertex $r$ to $u$ going through vertex $b$ or $d,$ this path splits the interior or exterior region of $C$ into 2 pieces. Note that it can split both, when the path crosses the curve $C$ several times, but we cannot discount such a situation to improve the bound.   
Hence we may use this path and a $T$-path from $r$ to one vertex from $\{u,v\}$ to build a Jordan curve $C'$ which allow us to save more than $\frac{1}{6}n-1$ vertices. See Figure~\ref{fig1} for the illustration.

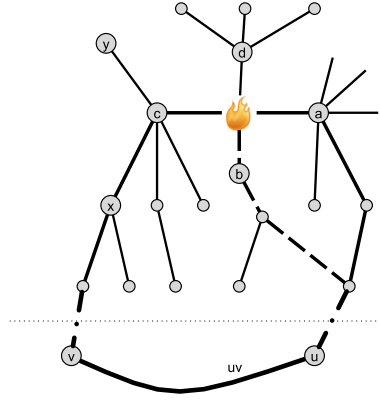
\begin{figure}[htbp]
\begin{minipage}[b]{0.5\textwidth}
\begin{flushright}
{\includegraphics[viewport=24 11 354 361,clip,width=2in]{Separator}}\\
\end{flushright}
\end{minipage}
\hspace{20pt}
\begin{minipage}[t]{0.38\textwidth}
\vspace{-110pt}
{\small
Solid edges are the edges of the spanning tree $T,$ bold if they belong to the curve $C.$ Dashed edges are the edges forming the  $ru$-path through $b$ (possibly but not necessarily belonging to the tree $T$). Note that the edges not belonging to $T$ are mostly not depicted.}
\end{minipage}
\caption{Illustration of the strategy.}\label{fig1}
\end{figure}

It may happen as well that in the graph $H$ there exists some shortest path from $r$ to $u$ going through vertex $c.$ From any of such paths let us choose this one which contains vertex $u_{i'}$ for the least possible $i',$ $1 < i' \le i.$
This path also splits either the interior or exterior region of $C$ into 2 pieces. If a piece not containing vertex $r$ has more than $\frac{n}{6}-1$ vertices we are done. Otherwise let $u' = u_{i'-1},$  $v' = u_{i'}$. Note that the curve $C'$ build from the $T$-path from $r$ to $u',$ edge $u'v'$ and the shortest $rv'$-path going through vertex $c$ has the property that $|\overline{\textrm{in} C'}| > \frac{1}{6}n+1$ and $|\overline{\textrm{ex} C'}| > \frac{1}{6}n+1.$ Moreover there is no shortest path from $r$ to $u'$ going through vertex $b,$ $c$ or $d.$ Hence we are ready to proceed to the final case with the curve $C',$ instead of $C.$ 

Suppose then that there is no such a ``shortcut'' and both closed interior as well as closed exterior of $C$ contains more than $\frac{n}{6}+1$ vertices. We have $\dist(b, u) \ge \dist(r, u), $ $\dist(c, u) \ge \dist(r, u)$ and $\dist(d, u) \ge \dist(r, u).$  Without loss of generality, assume that  $\deg(c) \le \deg(a)$ and that $b$ lies in the interior region of $C,$ while $d$ in the exterior. Note that terms ``interior'' and ``exterior'' depend on a particular embedding $H$ of $G.$ 

As $r$ has at most one neighbour of degree higher than 5, then one has $\deg(c) \le 5.$ Note that at most 3 neighbours of $c = v_1$ belong to $\overline{ \textrm{in} C}$ or at most 3 neighbours of $c$ belong to $\overline{\textrm{ex} C}.$ Suppose the exterior case holds (the interior one is analogous). One of these neighbours is $v_0 = r.$ Let the second (supposed) neighbour, say $x,$ belong either to $C$ ($x = v_2$ when $v_1 \neq v$) or to $\textrm{ex} C,$ while the 3rd one, if exists, is $y \in \textrm{ex} C.$ The firefighters' strategy is now the following. At their first round they protect $u_0 = a$ and $d.$ In the second one --- $x$ and $y$, then in the $t$th round ($t > 2$) --- $u_{(t-1)}$ and $v_{t}.$ This still allows to protect the exterior region of $C,$ because otherwise there would exist a path from $b$ or $c$ to $u$ of a length lower than $\dist(r, u).$ Hence, we are able to save more than $\frac{1}{6}n-1$ vertices ($r$ and $c$ are burned). See Figure~\ref{fig1} for details.
\end{proof}
Note that, in the above proof, considering the degree condition from Lemma~\ref{lem:2f} we do not count edges between neighbours of the vertex $r$. In particular, the following holds.
\begin{corollary}
Let $G$ and $r$ be as in Lemma~\ref{lem:2f}. Suppose that  $r$ is adjacent to 4 triangular faces. There is $sn_{2}(G, r) > \frac{1}{6}n -1$, provided that $r$ has at most one neighbour of degree higher than 7.
\end{corollary}

The rest of this section is devoted to calculate from the above lemmas the 2-surviving rate for a planar graph of a given average degree, to prove the bound given by Theorem~\ref{thm:main}. The bound is trivially true for graphs on not more than 4 vertices. Fix $n \in \N,$ $n \ge 5$ and $\epsilon > 0$ and let $G$ be any connected planar graph on $n$ vertices of average degree $\frac{9}{2} - \epsilon.$  

Let us partition the vertex set of $G$ into 4 subsets defined by the conditions:
\begin{eqnarray}
\nonumber X &=& \{v \in V(G) \colon \deg(v) \le 2\},\\
\nonumber Y &=& \left\{v \in V(G) \deg(v) = 3 \right\},\\
\nonumber Z &=& \left\{v \in V(G) \colon \deg(v) = 4 \wedge sn_{2}(G, v) \ge \frac{n-5}{6}\right\},\\
\nonumber W &=& V(G) \setminus (X \cup Y \cup Z).
\end{eqnarray}
By Observation~\ref{obs:simple}, one has $sn_{2}(G, v) \ge n-1$ for $v \in X$, $sn_{2}(G, v) \ge \frac{n-2}{3}$ for $v \in Y$, and $sn_{2}(G, v) \ge 2$ for $v \in W$. Therefore, we have that 
\begin{equation}\label{eq:sr2}
\rho_{2}(G) \ge \frac{1}{n^2} \bigg((n-1)|X| + \frac{n-2}{3}|Y| + \frac{n-5}{6}|Z| + 2|W|\bigg).
\end{equation}

Set $W$ contains all vertices of degree higher than 4. It may also contain vertices of degree 4, but (according to Lemma~\ref{lem:2f}) every such a vertex has at least 2 neighbours of degree at least 6. We claim that the average degree of vertices in the set $W$ is at least $\frac{9}{2}$. Indeed, suppose that the set $W_4 = \{v \in W \colon \deg(v) = 4\}$ is non-empty, let $W_6 = \{v \in W \colon \deg(v) \ge 6\}$, denote $q = |W_4|$, $t = |W_6|$, let $s = \sum_{v \in {W_6}} \deg(v)$. One has $t \le \frac{s}{6}$ and $2q \le s$ (each vertex in the set $W_4$ has at least 2 neighbours in $W_6$). The average degree of vertices in the set $W_4 \cup W_6$ is $$\frac{4q + s}{q+t} \ge \frac{4q+s}{q+\frac{s}{6}} = 4 + \frac{\frac{s}{3}}{q+\frac{s}{6}} \ge 4 + \frac{\frac{s}{3}}{\frac{s}{2}+\frac{s}{6}} = \frac{9}{2}.$$
Because vertices in the set $W \setminus (W_4 \cup W_6)$ have degree 5, the claim follows. 

It implies that 
\begin{equation}\label{eq:deg}
\Big(\frac{9}{2}-\eps\Big)n = \sum_{v \in V}\deg v \ge |X| + 3|Y| + 4|Z| + \frac{9}{2}|W|.
\end{equation}

One may then estimate inequality (\ref{eq:sr2}) reducing structural properties of graph $G$ to degree condition (\ref{eq:deg}) and solving the following linear program. %Note the correspondence between variables $x$, $u$, $z$, and $w$ and sizes of the sets $X$, $Y$, $Z$, and $W$.
\begin{eqnarray} \label{eq:lp}
\nonumber \textrm{minimize} & & \alpha = \frac{1}{n^2} \bigg((n-1)x + \frac{n-2}{3}y + \frac{n-5}{6}z + 2w\bigg)\\
 \textrm{with conditions} & & x+y+z+w = n,\\
\nonumber & & x + 3y + 4z + \frac{9}{2}w \le \Big(\frac{9}{2} - \eps \Big)n \\
\nonumber & & x, y, z, w \ge 0.
\end{eqnarray}

The solution is the following (see \cite{program} for calculations done for this paper, in particular for solutions of all linear programs).\\
\begin{footnotesize}
\noindent
\begin{tabular}{|r|l|c|c|c|c|c|}
\hline
& ~\hspace{1em}~\hspace{1em}~\hspace{1em}~\hspace{1em}~\hspace{1em}~\hspace{1em}~conditions & $x$ & $y$ & $z$ & $w$ & $\alpha$ \\
\hline
1. & $\eps \in (0, \frac{1}{2}],$ $n \le 17$ & 0 & 0 & $n$ & 0 & $\frac{1}{6}-\frac{5}{6n}$ \\
\hline
2. & $\eps \in (0, \frac{1}{2}],$ $17 \le n \le 35$ & 0 & 0 & $2\eps n$ & $(1-2\eps)n$ & $\frac{1}{3}\eps+\frac{1}{n}(2-\frac{17}{3}\eps)$  \\
\hline
3. & $\eps \in [\frac{1}{2}, \frac{3}{2}],$ $n \le 35$ & 0 & $(\eps-\frac{1}{2})n$ & $(\frac{3}{2}-\eps)n$ & 0 & $\frac{1}{6}\eps (1 + \frac{1}{n}) + \frac{1}{12} - \frac{11}{12n}$ \\
\hline
4. & $\eps \in (0, \frac{3}{2}],$ $n \ge 35$ & 0 & $\frac{2}{3}\eps n$ & 0 & $(1-\frac{2}{3}\eps)n$ & $\frac{2}{9}\eps+\frac{1}{n}(2-\frac{16}{9}\eps)$ \\
\hline
5. & $\eps \in [\frac{3}{2}, \frac{5}{2}]$ & $(\frac{1}{2}\eps-\frac{3}{4})n$ & $(\frac{7}{4}-\frac{1}{2}\eps)n$ & 0 & 0 & $\frac{1}{3}\eps-\frac{1}{6}(1+\frac{\eps}{n})-\frac{5}{12n}$ \\
\hline 
\end{tabular}
\end{footnotesize}

In each case $\alpha \ge \frac{2}{9}\eps - \frac{1}{n}$. Moreover, in case 2 and case 4 (for $\eps \in (0,\frac{9}{8}]$) one has $\alpha \ge \frac{2}{9}\eps$. For the first case note that $\rho_2(G)$ can be easily bounded directly: protecting 2 vertices in the first round firefighters protect at least $\frac{2}{17} > \frac{2}{9}\eps$ vertices of the graph, no matter where the fire starts. For the case 3 and for $\eps \in (\frac{1}{2}, 1]$ at first note that $\frac{1}{6}\eps (1 + \frac{1}{n}) + \frac{1}{12} - \frac{11}{12n} \ge \frac{2}{9}\eps$ for $n \ge 27$. Then, notice that protecting 2 vertices at the first step it is enough to obtain the bound $\rho_2(G) \ge \frac{2}{9}\eps$ for $n \le 9$.

Thus, only the bound for $\eps \in (\frac{1}{2}, 1]$ and $n$ between 10 and 26 remains to be proven.  Notice that for small graphs Lemma 3.2 is  no longer useful. Indeed, when the fire breaks out at vertex $r$ such that $\deg(r) \le n-3$, firefighters can protect 4 vertices: 2 in the first round and 2 in the second one, because at least 2 vertices are not neighbouring $r$. One can then modify the definition of the set $Z$ to $$Z' = \{v \in V(G) \colon sn_2(v) \ge 4\} \setminus (X \cup Y).$$ Then the set $W'$ contains vertices of degree at least $n-2$, which leads to estimate the bound for $\rho_2(G)$ with the following linear program.
\begin{eqnarray} \label{eq:lp3}
\nonumber \textrm{min} & & \alpha = \frac{1}{n^2}\Big((n-1)x + \frac{n-2}{3}y+4z' + 2w'\Big)\\
 \textrm{with conditions} & & x+y+z'+w' = n,\\
\nonumber & & x + 3y + 4z' + (n-2)w' \le \Big(\frac{9}{2} - \eps \Big)n \\
\nonumber & & x, y, z', w' \ge 0.
\end{eqnarray}
For $15 \le n \le 26$ it has the solution $\alpha = \frac{1}{3} \eps (1 -\frac{14}{n})-\frac{1}{6}+\frac{19}{3n}$; for $10 \le n \le 14$ --- the solution $\alpha = \frac{n+\eps}{3(n-5)}-\frac{17n+16\eps-44}{6n(n-5)}.$ Both solutions, for appropriate values of $n$, satisfy the inequality $\alpha \ge \frac{2}{9}\eps,$ which closes the proof of Theorem~\ref{thm:main}.

\section{Remarks} \label{sec:remarks}
Noting the results on the 2-surviving rate for planar graphs without 3-, 4-, 5- or 6-cycles, respectively~\cite{esperet,kwzhang,kwwu,fww}, one may ask the following question:
\begin{problem}
Does it exist a function $f \colon \N\setminus\{1,2\} \to \mathbb{R}^+$ such that for any $k \ge 3$ and any non-trivial planar graph $G$ without $k$-cycles  there is $\rho_2(G) > f(k)?$
\end{problem}
Of course, such a function does exist provided there is a strategy for the positive 2-surviving rate on any planar graph, which is still open. 

\section{Acknowledgement} \label{sec:ack}
As written in the abstract, this paper is the corrected version of~\cite{pg_planarfire2}. There was some technical/calculation error in the proof of Theorem~~\ref{thm:main}, namely the solution of an auxiliary parameterised linear program (analogous to (\ref{eq:lp})) was partially incorrect which caused a gap in the proof for graphs on less than 36 vertices. The corrected proof is presented above while Theorem~\ref{thm:main} remains unchanged. Also Corollary~\ref{cor:4cycles} required some modification: in the original paper the bound for average degree for planar graphs without 4-cycles was miscalculated.

The author would like to thank Bartosz Walczak for his careful reading of~\cite{pg_planarfire2} and pointing out the errors. 

\end{document}